\documentclass{amsart}
\usepackage{latexsym,amsmath,amsopn,amssymb,amsthm,amsfonts}

\newtheorem*{theorem*}{Theorem}
\newtheorem{theorem}{Theorem}

\newtheorem{lemma}{Lemma}
\newtheorem*{lemma*}{Lemma}
\newtheorem*{conj*}{Conjecture}

\theoremstyle{definition}
\newtheorem*{example*}{Example}

\theoremstyle{remark}

\newcommand{\RR}{\mathbb{R}}
\newcommand{\CC}{\mathbb{C}}
\newcommand{\CP}{\mathbb{C}P}

\newcommand{\ZZ}{\mathbb{Z}}

\newcommand{\vv}[2]{\begin{pmatrix}#1\\#2\end{pmatrix}}
\newcommand{\vb}[2]{\begin{bmatrix}#1\\#2\end{bmatrix}}
\newcommand{\matriza}[4]{\begin{pmatrix}#1&#2\\#3&#4\end{pmatrix}}
\newcommand{\case}[1]{({\sc #1})}

\newcommand{\im}{\mathrm{Im}\,}
\newcommand{\de}{\partial}
\newcommand{\chp}[2]{\frac{\partial#1}{\partial#2}}

\begin{document}

\title{Theta-functions on $T^2$-bundles over $T^2$ with the Euler
class zero}
\author{Dmitry V. Egorov}
\address{
Ammosov Northeastern federal university, \newline Kulakovskogo str. 48, 677000, Yakutsk, Russia}%
\email{egorov.dima@gmail.com}%
\date{}

\thanks{This work was supported in part by Russian Foundation for Basic
Research (06-01-00094-a).}

\begin{abstract}
We construct an analogue of the classical theta-function on an
Abelian variety for closed $4$-dimensional symplectic manifolds
which are $T^2$-bundles over $T^2$ with the zero Euler class. We use
our theta-functions for a canonical symplectic embedding of these
manifolds into complex projective spaces (an analogue of the
Lefschetz theorem).
\end{abstract}

\maketitle

\section{Introduction}

We proposed in \cite{egorov1} a construction of theta functions on
the Kodaira–-Thurston manifold. Here we present a different approach
to their construction on closed $4$-dimensional symplectic
manifolds, namely $T^2$-bundles over $T^2$ with the zero Euler
class. From the geometric viewpoint, the classical theta function on
an Abelian variety is a section of a holomorphic line bundle over
the complex torus. The Lefshetz theorem states that a section of a
suffciently high tensor power of the bundle determines a
complex--analytic embedding of the Abelian variety into some complex
projective space. We aim at generalizing this construction to the
bundles whose fibers and bases are $1$-dimensional complex tori,
while the Euler class is zero. We introduce some analogues of the
classical theta-functions as sections of complex line bundles over
these bundles. We have to relinquish the holomorphic embedding since
these bundles, generally speaking, lack not only a K\"{a}hler
structure, but even a complex structure. Nevertheless, they are
symplectic manifolds. We will construct theta functions so that a
symplectic analogue of the Lefschetz theorem holds for them: theta
functions with characteristics, which are sections of tensor powers
of the defining line bundle, determine a symplectic embedding of the
manifold into $\CP^k$ (for sufficiently high tensor powers).

The Kodaira–-Thurston manifold can be viewed as a $T^2$-bundle over
$T^2$ in two distinct ways (the bundles are not isomorphic), one of
which we consider in the article. Therefore, the theta functions on
the Kodaira–-Thurston manifold introduced in this article differ by
construction from the theta functions of \cite{egorov1}, and it
turns out that they coincide with the theta functions of Kirwin and
Uribe \cite{KU}, whose approach is representation--theoretic. We
discuss this coincidence in more detail in Subsection $4.4$. In
Section $2$ we recall the necessary facts of the classical theory of
theta-functions and in Section $3$ we describe the bundles with
which we are to work. In Section $4$ we define theta-functions on
bundles and study some of their properties, in Section $5$ we
construct an embedding of those bundles into complex projective
spaces (Theorem $1$), and in Section $6$ we prove that the embedding
is symplectic (Theorem $2$). It would be interesting to find out how
far the analogy with the classical theta-functions goes; for
instance, whether the constructed theta-functions are related to
number theory (see \cite{Mumford} for instance) or nonlinear
equations and secant formulas (see the survey \cite{T}). The author
is grateful to I. A. Taimanov for stating the problem and A. E.
Mironov for useful discussions.

\section{The classical theta-function}

Recall some useful facts concerning the theta-function on the
$1$-dimensional complex torus.

Consider the formal series
$$
\theta(z,\tau) =
e^{\pi iz + {\pi i\tau}/{4} + {\pi i}/{2}}\sum_{k\in\ZZ}{e^{2\pi i kz + \pi ik(k+1)\tau + \pi ik}}.
$$
For $\mathrm{Im}\,\tau>0$ this series converges in every compact
domain in $\mathbb{C}$ and defines an entire function. With respect
to translations the theta function behaves as
\begin{eqnarray}
\label{classperiod}
\theta(z+1, \tau)& = &-\theta(z, \tau), \\
\label{classperiod2}
\theta(z+\tau, \tau)& = & -e^{-2\pi iz - \pi i\tau}\theta(z,\tau).
\end{eqnarray}

Our choice of this particular theta function, which is
$\theta_{11}(z, \tau)$ in the notation of Mumford, is related to the
fact that it becomes multiplied by an exponential under all modular
transformations:
\begin{equation}
\label{modular}
\theta(\frac{z}{c\tau + d}, \frac{a\tau+b}{c\tau + d}) = \zeta(c\tau+d)^{1/2}\exp(\frac{cz^2}{c\tau+d})\cdot\theta(z,\tau), \quad ad-bc = 1.
\end{equation}
Here the nonzero constant $\zeta$ depends on the matrix
$\begin{pmatrix}
                                                            a & b \\
                                                            c & d \\
                                                          \end{pmatrix}$.

A generalization of the concept of theta-function is the concept of
theta-function of degree $k$, where $k$ is some positive integer. A
theta-function of degree $k$ is an entire function on $\CC$ with the
periodicity properties
$$
\theta_k(z+1, \tau) =  (- 1)^k\theta_k(z, \tau),
$$
$$
\theta_k(z+\tau, \tau) =  (-1)^k e^{(-2\pi iz - \pi i \tau )k}\:\theta_k(z,\tau).
$$
It is not difficult to verify that theta-functions of degree $k$
span a linear space of dimension $k$ denoted by $\mathcal{L}_k$.

Multiplying theta functions, we obtain a theta-function of a higher
degree. Take a tuple $\{\alpha_i\}^k_{i=1}$ of constants whose sum
is equal to zero. Then
$$
\prod_{i=1}^k{\theta(z+\alpha_i,\tau)} \in \mathcal{L}_k
$$

The theta function is equal to zero at the point $z = 0$ modulo the
lattice $\ZZ + \tau \ZZ$. The fundamental domain of the lattice
contains a unique zero counting the multiplicity.

The theta function satisfies the following partial differential
equation:
\begin{equation}
\label{classheat} \chp{\theta(z,\tau)}{\tau} = \frac{1}{4\pi
i}\frac{\de^2\theta(z,\tau)}{\de z^2}.
\end{equation}

\section{Bundles}

The manifolds that are the total spaces of the bundles whose fiber
and base are $2$-dimensional tori are classified in \cite{Fukuhara}.
Let us describe the bundles with the zero Euler class.

Take two commuting matrices $A$ and $B$ in $SL(2,\ZZ)$. Denote by
$\vb{s}{t}$ the point of $T^2=\RR^2/\ZZ^2$ corresponding to
$\vv{s}{t}\in\RR^2$. Then the bundle corresponding to
$\left\{A,B\right\}$ arises upon taking the quotient of
$T^2\times\RR^2/\sim$ by the action of $A$ and $B$:
$$
\left(\vb{s}{t},\vv{x+1}{y}\right)\sim
\left(\left[A\vv{s}{t}\right],\vv{x}{y}\right)
$$
and
$$
\left(\vb{s}{t},\vv{x}{y+1}\right)\sim
\left(\left[B\vv{s}{t}\right],\vv{x}{y}\right).
$$
A complete list of these bundles, up to a diffeomorphism of the
total space, is shown in the table, where I is the identity matrix.

\begin{center}
\begin{tabular}{|c|l|l|}
  \hline
  & $\left\{A,B\right\}$ & Notes \\
  \hline\hline
  ({\sc A}) & $\left\{I, I\right\}$ =  $T^4$ & $4$-torus \\
  \hline
  ({\sc B}) &(1) $\left\{\matriza{0}{-1}{1}{-1}, I\right\}$ & Hyperelliptic surface\\
   &   (2)  $\left\{\matriza{0}{-1}{1}{0}, I\right\}$ & \\
   &   (3)  $\left\{\matriza{1}{-1}{1}{0}, I\right\}$ & \\
   &   (4)  $\left\{-I, I\right\}$ &  \\
   \hline
  ({\sc C}) & $\left\{\matriza{1}{k}{0}{1}, I\right\}$, $k\neq 0$ & Kodaira--Thurston manifold\\
  \hline
  ({\sc D}) & $\left\{\matriza{-1}{k}{0}{-1}, I\right\}$, $k\neq 0$ & \\
  \hline
  ({\sc E}) & $\left\{\matriza{1}{k}{0}{1}, -I\right\}$, $k\neq 0$  &  \\
   \hline
  ({\sc F}) & $\left\{A, I\right\}$, |$\mathrm{tr}\, A| > 2$ &  \\
  \hline
  ({\sc G}) & $\left\{A, -I\right\}$, $\mathrm{tr}\, A > 2$ &  \\
  \hline
\end{tabular}
\end{center}
A survey devoted to the manifolds of type ({\sc B}) can be found in
\cite{Hasegawa} for instance.

Verify that all these bundles are homogeneous spaces of real Lie
groups. Take a real Lie group G with coordinates $(x, y, s, t)$ and
the linear representation
$$
\rho: (x,y,s,t)\to \left(
\begin{tabular}{c|c}
  $A^xB^y$ & $\begin{array}{c}
                s \\
                t
              \end{array}
  $ \\
  \hline
  $\begin{array}{cc}
    0 & 0
  \end{array}$
   & 1 \\
\end{tabular}
\right).
$$
Denote by $\Gamma$ the discrete subgroup of elements with integer
coordinates. Then the homogeneous space $\Gamma\backslash G$ is the
total space of the bundle corresponding to $\left\{A,B\right\}$.
Observe that $\rho$ is a faithful representation not of G, but
rather of $\Gamma\backslash G$.

It is obvious that we can proceed to define these bundles as the
quotient manifolds of $\RR^4$ by the action of the group $\Gamma$
with the following generators (when $B = I$):
\begin{eqnarray}
\label{gammaA}
a&: & (x+1, y, \alpha s + \beta t, \gamma s + \delta t),\\
\label{gammaB}b&: & (x, y+1, s, t), \\
\label{gammaC}c&: & (x, y, s+1, t), \\
\label{gammaD}d&: & (x, y, s, t+1),
\end{eqnarray}
where $ A = \matriza{\alpha}{\beta}{\gamma}{\delta}. $ If $B=-I$,
then
\begin{eqnarray}
\label{gammaBB}b&:& (x, y+1, -s, -t).
\end{eqnarray}

Knowing a faithful representation $\Gamma\backslash G$, we can
calculate the generators of the algebra of left-invariant forms. As
usual
$$ \rho^{-1}d\rho = \left(
\begin{tabular}{c|c}
  $C dx + Ddy$ & $A^{-x}B^{-y}\vv{ds}{dt}$ \\
  \hline
  $\begin{array}{cc}
    0 & 0
  \end{array}$
   & 0 \\
\end{tabular}
\right),
$$
where $C$ and $D$ are some constant matrices. Therefore, $1$-forms
$dx$, $dy$, $\omega_1$, $\omega_2$, where
$$
\vv{\omega_1}{\omega_2} = A^{-x}B^{-y}\vv{ds}{dt},
$$
are the generators.

\section{Theta-functions on bundles}

\subsection{The definition of theta-function}
Henceforth we denote by $M$ the total space of bundles whose the
fiber and base are $2$-dimensional tori and whose Euler class is
zero.

Introduce a function on the universal covering $R^4$ of $M$
formally:
\begin{equation}
\label{my_theta} {\vartheta_M}(x,y,s,t) = \theta(s+\omega
t,\omega)\cdot \theta(x+iy,i).
\end{equation}
Here $\theta(z,\tau)$ is the classical theta function of the
variable $z$ with period $\tau$ (see Section 2). Specify the
function $\omega$ depending on the bundle (writing $i = \sqrt{-1}$):
%~\\
\begin{center}
\begin{tabular}{c|l}
   & $\omega$  \\
  \hline\hline
  \case{B1}, \case{B3} & $(-1 + \sqrt{-3})/2$ \\
  \case{B2}, \case{B4} & $i$ \\
  \case{C}, \case{E} & $-kx +i$ \\
  \case{D} & $kx +i$ \\
  \case{F}, \case{G}  & $(\lambda^{-x}v^+ + i\lambda^{x}v^-)/(\lambda^{-x}u^+ + i\lambda^{x}u^-)$
\end{tabular}
\end{center}
%~\\
Here $\lambda$ and $\lambda^{-1}$ are the eigenvalues of $A^T$,
while $(u^+,v^+)^T$ and $(u^- ,v^-)^T$ are eigenvectors of the
transposed matrix $A$:
$$ A^T
\begin{pmatrix}
u^+ \\
v^+
\end{pmatrix}
= \lambda \begin{pmatrix}
u^+ \\
v^+
\end{pmatrix}
,\quad A^T \begin{pmatrix}
u^- \\
v^-
\end{pmatrix}
 = \lambda^{-1} \begin{pmatrix}
u^- \\
v^-
\end{pmatrix}.
$$
the vectors are normalized so that $u^+v^- - u^-v^+ = 1$.

\begin{lemma}
$\im \omega >0$ for all bundles.
\end{lemma}
{\sc Proof.}\quad Only the case of the last function is not obvious:
$$
\im \omega = \frac{u^+v^- - u^-v^+}{(\lambda^{-x}u^+)^2 +
(\lambda^{x}u^-)^2} >0.
$$
The proof of the lemma is complete.

By Lemma $1$ the series defining the theta function converges, and
the function $\theta_M$ is well-defined.

\begin{lemma}
Under the action of the generator \eqref{gammaA} of $\Gamma$:
$$
(x, y, s, t)\to(x+1, y, \alpha s +\beta t, \gamma s +\delta t)
$$
functions $\omega$ and $s+\omega t$ transform as:
$$
\omega \rightarrow \frac{\alpha\omega - \beta}{-\gamma\omega +
\delta },\quad s+\omega t \rightarrow \frac{s+\omega t}{-\gamma
\omega +\delta}.
$$
\end{lemma}
{\sc Proof.}\quad Direct calculations.

Take a basis $\{\theta_k^p(z,\tau)\}_{p=1}^k$ for the space of
classical theta-functions of degree $k$ (see Section $2$). Define
the space of theta-functions of degree $k$ on the manifold $M$ as
the linear span of pairwise products of ordinary basis
theta-functions of degree $k$ on the fiber and base:
$$
\theta^p_k(s+\omega t,\omega)\cdot\theta_k^q(x+iy,i),\quad p,q =
1\ldots k
$$
Denote this space by $\mathcal{L}_k$. Observe that the dimension of
$\mathcal{L}_k$ is equal to  $k^2$. The theta function of degree $1$
is precisely ${\vartheta_M}$.

\subsection{${\vartheta_M}$ is a section of a complex line bundle}
In order to show that the function $\theta_M$ is a section of a
complex line bundle over $M$, recall that the sections are in a
bijective correspondence with the functions $f$ on the universal
covering such that $f(\lambda\cdot u) = e_\lambda(u)f(u)$, where
$\lambda$ is an element of the lattice $\Gamma$ and  $e_\lambda(u)$
are the multipliers; i.e., nonzero functions $e_\lambda : \RR^4
\rightarrow \CC^{*} $ such that
$$
e_{\lambda}(\mu\cdot u)\:e_{\mu}(u) = e_{\lambda\mu}(u),\quad
\lambda,\mu\in \Gamma
$$
$$
e_{0}(u) = 1.
$$
The multipliers determine a complex line bundle over $M$ such that
the direct product $\RR^4\times\CC$ can be modded out by the action
of $\Gamma$:
$$
(u,w) \sim (\lambda\cdot u,
e_{\lambda}(u)w),\quad u\in\RR^4,w\in\CC,\lambda \in\Gamma.
$$
Consider the behavior of ${\vartheta_M}$ under the action of the
generators \eqref{gammaA}-\eqref{gammaD} of $\Gamma$ (with $B=I$):
\begin{eqnarray}
&&\label{mult1}{\vartheta_M}(x+1, y, \alpha s +\beta t, \gamma s +\delta t) =\\
&&\nonumber\qquad\qquad\qquad \zeta(-\gamma \omega + \delta)^{\frac{1}{2}}\exp(\frac{-\gamma (s+\omega t)^2}{-\gamma \omega + \delta})\cdot {\vartheta_M}(x,y,s,t),\\
&&\label{mult2}{\vartheta_M}(x, y+1,  s, t) =  -\exp(-2\pi i (x+iy)+\pi )\cdot{\vartheta_M}(x, y, s, t),\\
&&\label{mult3}{\vartheta_M}(x, y, s+1, t) =  -{\vartheta_M}(x, y, s, t),\\
&&\label{mult4}{\vartheta_M}(x, y, s, t+1) = -\exp(-2\pi i (s+\omega
t)-\pi i\omega )\cdot {\vartheta_M}(x, y, s, t).
\end{eqnarray}
In the case $B=-I$ we should replace  \eqref{mult2} by
\begin{equation}
\quad \label{mult5}{\vartheta_M}(x, y+1,  -s, -t) =  \exp(-2\pi i
(x+iy)+\pi )\cdot{\vartheta_M}(x, y, s, t).
\end{equation}
We used the periodicity properties
\eqref{classperiod}-\eqref{modular} of the classical theta function
and Lemma 2. These formulas imply that $\vartheta_M$ is a section of
the bundle determined by the multipliers \eqref{mult1} -
\eqref{mult5}.

In order to verify that the construction of the bundle is
well-defined, given these multipliers, we must verify that the
nontrivial relations \eqref{gammaA}-\eqref{gammaBB}:
$$
[a,c] = c^{1-\delta}d^\gamma,\quad [a,d] = c^\beta
d^{1-\alpha},\quad [g,h] = g^{-1}h^{-1}gh,
$$
between the generators \eqref{gammaA}-\eqref{gammaBB} of $\Gamma$
imply identities on the multipliers. This is obvious considering
that all multipliers are determined by the behavior of the same
function.

\subsection{A multiplicative property of ${\vartheta_M}$}
Introduce the action of $\zeta = (\lambda,\mu)\in \CC^2$ on
${\vartheta_M}$ as
\begin{equation}
\label{preprod} (\zeta\cdot {\vartheta_M})(x,y,s,t) =
\theta(s+\omega t + \lambda,\omega)\cdot\theta(x+iy + \mu,i).
\end{equation}
Take a tuple $\zeta_i = (\lambda_i,\mu_i),\ i = 1,\ldots , k$ of
constant vectors in $\CC^2$ whose sum is equal to zero.  As for the
classical theta function, it is desirable that the product
\begin{equation}
\label{product} \prod_{i=1}^k{(\zeta_i\cdot{\vartheta_M})(x,y,s,t)}
\end{equation}
be a theta function of degree  $k$. This property of the theta
function is key for the proof of our theorem on the embedding into a
complex projective space. It is easy to verify that the multiplier
\eqref{mult1} of the form $\exp(-\gamma(s+\omega t)^2)$ is an
obstruction. If this multiplier is nontrivial, i.e.,  $\gamma\neq
0$, then we require also that
\begin{equation}
\label{postprod} \sum_{i=1}^k{\lambda_i^2} =0.
\end{equation}
It is not difficult to verify then that the product  \eqref{product}
is a theta function of degree $k$.

\subsection{Relation to the theta functions of Kirwin and Uribe}
\label{sectionKU} Let us recall how theta functions on the
Kodaira–-Thurston manifold were introduced in  \cite{KU}. Given a
square-integrable function $f:\RR^2\to\CC$ for all $m,n =
0,1,\ldots, 2k-1$ the functions $\vartheta_k^{m,n} f:\RR^4\to \CC$
are defined as
\begin{multline*} (\vartheta_k^{m,n}f)(x,y,z,t) =\\
e^{-2\pi i [my - n(z+xy)]-4\pi i k zx}\sum_{a,b\in\ZZ}e^{2\pi i
nya-4\pi i k(by-za-y(x+a)^2/2)}f(x+a,t+b).
\end{multline*}
These functions satisfy the pseudo-periodicity
conditions\begin{eqnarray*}
(\vartheta_k^{m,n}f)(x+1,y,z,t) &=& (\vartheta_k^{m,n}f)(x,y,z,t),\\
(\vartheta_k^{m,n}f)(x,y+1,z-x,t) &=& e^{-2\pi i kx^2}(\vartheta_k^{m,n}f)(x,y,z,t),\\
(\vartheta_k^{m,n}f)(x,y,z+1,t) &=& e^{4\pi i kx}(\vartheta_k^{m,n}f)(x,y,z,t),\\
(\vartheta_k^{m,n}f)(x,y,z,t+1) &=& e^{4\pi i
ky}(\vartheta_k^{m,n}f)(x,y,z,t).
\end{eqnarray*}

Verify that for a certain choice $f:\RR^2\to\CC$ the function
$\vartheta_k^{m,n}f$ turns out with some degree of precision to be
the theta function on the Kodaira-–Thurston manifold defined in this
article.

Take $k=1$, $m=0$, $n=0$, $f = g(t)\cdot h(x)$:
$$
g(t) = e^{-2\pi t^2}, \quad h(x) =  e^{-2\pi x^2}.
$$
Then
\begin{multline}
\label{KUtheta} (\vartheta_k^{m,n}f)(x,y,z,t) =\\ e^{-4\pi i zx-2\pi
t^2 -2\pi x^2}\cdot \theta(2(z+(y+i)x), 2(y+i))\cdot
\theta(2(-y+it), 2i),
\end{multline}
where $\theta(z,\tau)$ is the classical theta function with
characteristics $[0,0]$.

Recall that the theta function on the Kodaira–-Thurston manifold
defined here in  \eqref{my_theta} is of the form
$$
\theta[\frac{1}{2},\frac{1}{2}](s+\omega t,\omega)\cdot
\theta[\frac{1}{2},\frac{1}{2}](x+iy,i),\quad \omega = -kx+i,\quad
k\in\ZZ \backslash\{0\},
$$
where $\theta[a,b](z,\tau)$ is the classical theta function with
characteristics  $[a,b]$. We chose a somewhat different theta
function because of its invariance under modular transformations.
Kirwin and Uribe multiply the argument and the period by 2 for
seemingly the same purpose. The change of variables
$$
x' = t,\: y' = -x,\: z' = s,\: t' = y,\: k=1,
$$
carries the function
$$
\theta[\frac{1}{2},\frac{1}{2}](s+\omega t,\omega)\cdot
\theta[\frac{1}{2},\frac{1}{2}](x+iy,i),\quad \omega = -kx+i
$$
into
$$
\theta[\frac{1}{2},\frac{1}{2}](z'+(y'+i)x', y'+i)\cdot
\theta[\frac{1}{2},\frac{1}{2}](-y'+it', i),
$$
which coincides with \eqref{KUtheta} up to the translations of
characteristics and multiplication by an exponential.

\section{Embedding into a Complex Projective Space}
Enumerate the basis theta functions $\{\sigma_i\}_{i=1}^{k^2}$ in
$\mathcal{L}_k$. Then
$$
\varphi_k = \left(\sigma_1 , \sigma_2 , \ldots , \sigma_{k^2}\right)
$$
is a well-defined mapping of the manifold $M$ into $\CP^{k^2-1}$.

Recall  that we denote the entries of the monodromy matrix $A$ by $
\matriza{\alpha}{\beta}{\gamma}{\delta}$, and   for $\gamma\neq 0$
we require the fulfilment of the additional condition
\eqref{postprod}.

\begin{theorem}
The mapping $\varphi_k$ is an embedding provided that
\begin{itemize}
\item[(a)] $k\geq 4$ when $\gamma\neq 0$;
\item[(b)]$k\geq 3$ when $\gamma=0$.
\end{itemize}
\end{theorem}

{\sc Proof.}\quad Let us prove claim  (a) of the theorem for $k=4$.
It should be clear from the proof how to deal with the remaining
cases.

First establish the injectivity of $\varphi_k$. We will follow the
proof of the classical Lefschetz embedding theorem for Abelian
varieties (see the exposition of it in \cite[Chapter 2, Theorem
1.3]{Mumford}).

Observe that the space of theta functions of degree $k$ consists of
the global sections of the $k$th tensor power of the bundle
determined by the multipliers \eqref{mult1}-\eqref{mult5}.

If it is true that for all points $u\neq v\in M$ there exists a
section $\sigma\in \mathcal{L}_k$ such that $\sigma(u)=0$ and
$\sigma(v)\neq 0$, then $\varphi_k$ is injective. Indeed, suppose
that the mapping "glues" together $u$ è $v$. Since $\varphi_k$ is
made up from the basis sections of $\mathcal{L}_k$, it follows that
$\sigma(v) = \zeta \cdot \sigma(u)$ for every section
$\sigma\in\mathcal{L}_k$, where $\zeta$ is some non-zero constant.
If $\sigma$ is a section    satisfying the condition indicated above
then we arrive at a contradiction. Observe also that, when this
condition is fulfilled at every point $u\in M$ not all sections
vanish at $u$.

We will seek a theta function of degree $4$ as a product $\sigma =
f\cdot g$ of two functions:
\begin{multline}
\label{section} f(s+\omega(x) t,x,\alpha,\beta) = \theta(s+\omega
t+\alpha,\omega)\theta(s+\omega t+\beta,\omega)\times\\\times
\theta(s+\omega t+\gamma,\omega)\theta(s+\omega t+\delta,\omega),
\end{multline}
\begin{multline}
\label{section2} g(x+iy,\alpha',\beta',\gamma') =  \theta(x+iy
+\alpha',i)\theta(x+iy+\beta',i)\times\\\times\theta(x+iy+\gamma',i)\theta(x+iy+\delta',i).
\end{multline}Here
\begin{equation}
\gamma = \frac{1}{4}\left({-2(\alpha+\beta) +
\sqrt{-4(\alpha+\beta)^2 - \alpha^2-\beta^2}}\right),
\end{equation}
\begin{equation}
\delta = \frac{1}{4}\left({-2(\alpha+\beta) -
\sqrt{-4(\alpha+\beta)^2 - \alpha^2-\beta^2}}\right),
\end{equation}
\begin{equation}
\label{sectionend} \delta' = -\alpha' - \beta'-\gamma'.
\end{equation}
In Subsection $4.3$  we showed that $f\cdot g$ is indeed a
theta-function of degree $4$ on $M$.

Denote the coordinates of  $u$ and $v$ by $(x,y,s,t)$ and
$(x',y',s',t')$ respectively. Choose $\alpha'$ so that
$\theta(x+iy+\alpha') = 0$. Now choose $\beta',\gamma',\delta'$ so
that the remaining factor in the definition of $g$ does not vanish
at $v$:
$$
\theta(x'+iy'+\beta')\theta(x'+iy'+\gamma')\theta(x'+iy'+\delta')\neq
0.
$$
We can achieve this since the zeros of the theta function are
isolated. By a small perturbation of $\alpha,\beta,\gamma,\delta$ we
can also achieve the non-vanishing of $f$ at $v$.

Observe that now it is easy to explain why in claim  (a) $k$ must be
at least $4$. For $k=3$, when there is no  $\delta'$, the constants
$\beta'$,$\gamma'$ are functions of $\alpha'$ because of
\eqref{postprod}, and it is impossible to achieve the non-vanishing
indicated above.

The constructed section would solve the problem provided that
$\theta(x'+iy'+\alpha')\neq 0$. Suppose that it is so. Since the
classical theta function has a unique zero in the fundamental domain
of the lattice formed by its periods, it follows that $x = x', y =
y'$. The equality is understood modulo the lattice, but without
restricting the generality we may assume that $u$ and $v$ lie in the
fundamental domain; i.e., the unit cube $0\leq x,y,s,t <1$.

Choose $\alpha$ so that  $\theta(s +\omega(x) t + \alpha,\omega(x))
= 0$. Observe that then $\theta(s' +\omega(x') t' +
\alpha,\omega(x')) \neq 0$, for otherwise $u=v$. Pick
$\beta,\gamma,\delta$ so that  $f(v)\neq 0$, while
$\alpha',\beta',\gamma',\delta'$ so that  $g(v)\neq 0$.

Therefore, we have constructed a required section and proved the
injectivity of $\varphi_k$.

Let us now prove that the rank of $\varphi_k$ is maximal. We will
follow the proof of the Lefschetz theorem in \cite{T}. To start off
show that the rank of the mapping is maximal provided that the rank
(over $\CC$) of the following matrix is maximal
\begin{equation*}
J =
\begin{pmatrix}
  \sigma_1 & & \ldots & & \sigma_{k^2} \\
  \partial_x\sigma_1 & & \ldots & & \partial_x\sigma_{k^2} \\
  \partial_y\sigma_1 & & \ldots & & \partial_y\sigma_{k^2} \\
  \partial_s\sigma_1 & & \ldots & & \partial_s\sigma_{k^2} \\
  \partial_t\sigma_1 & & \ldots & & \partial_t\sigma_{k^2} \\
\end{pmatrix}.
\end{equation*}

Observe that  $\varphi_k$, expressed in the homogeneous coordinates,
is the composition of a mapping $\widetilde{\varphi}_k$ into
$\mathbb{C}^{k^2}$ and the subsequent projection $\pi:
\mathbb{C}^{k^2}\backslash \{0\}\rightarrow\mathbb{CP}^{k^2-1}$. It
is obvious that the differential of $\widetilde{\varphi}_k$
coincides with the submatrix of $J$ resulting upon the removal of
the first row.

Suppose now that at a point $u^*\in {M}$ the first row of $J$ is a
linear combination of the remaining rows. This means that the
radius-vector $\widetilde{\varphi}_k(u^*)$ is collinear to the image
of some tangent vector at $u^*$. Since $\pi$ projects along the
complex lines passing through the origin, it follows that the kernel
of the differential of $\pi$ consists of precisely those vectors.
Consequently, the maximality of the rank of $J$ is a necessary and
sufficient condition for the maximality of the rank of $\varphi_k$.

Transform the matrix $J$ into a form convenient for us. The rank of
the matrix
\begin{equation*}
\widetilde{J} =
\begin{pmatrix}
  \sigma_1 & & \ldots & & \sigma_{k^2} \\
  (\partial_x-i\partial_y)\sigma_1 & & \ldots & & (\partial_x-i\partial_y)\sigma_{k^2} \\
  (\partial_s+\frac{\partial_t}{\omega})\sigma_1 & & \ldots & & (\partial_s+\frac{\partial_t}{\omega})\sigma_{k^2} \\
  (\partial_x+i\partial_y)\sigma_1 & & \ldots & & (\partial_x+i\partial_y)\sigma_{k^2} \\
  (\partial_s-\frac{\partial_t}{\omega})\sigma_1 & & \ldots & & (\partial_s-\frac{\partial_t}{\omega})\sigma_{k^2} \\
\end{pmatrix}.
\end{equation*}
coincides with the rank of $J$. Recall that for a holomorphic
function $f$ of a complex variable $w = u + iv$ we have
$$
\chp{f}{w} = \frac{1}{2}(\de_u-i\de_v)f,\quad \chp{f}{\bar{w}} =
\frac{1}{2}(\de_u+i\de_v)f = 0.
$$
The last two rows of $\widetilde{J}$ appear in the Cauchy–-Riemann
conditions. Since for every $x$ the sections $\sigma_j$ expand into
series  $s+\omega(x)t$, the last row of  $\widetilde{J}$ is always
zero.

Observe that if $\omega=\mathrm{const}$ (which we have for the theta
functions on the bundles of type \case{B}) then the sections are
holomorphic, and we can use the proof of the classical Lefschetz
theorem.  We will assume that   $(\de_x+i\de_y)\theta_M\neq 0$.

Assume that the rank of $\widetilde{J}$ (over $\CC$) at some fixed
point $u^*=(x^*,y^*,s^*,t^*)\in M$ is less than $4$. This means that
there exists a nontrivial tuple $a,b,c,d$ such that
\begin{multline*}
a\sigma_j(u^*) + \frac{b}{2}(\partial_x-i\partial_y)\sigma_j(u^*) +
\frac{c}{2}(\partial_s+\frac{\partial_t}{\omega})\sigma_j(u^*) +
\frac{d}{2}(\partial_x+i\partial_y)\sigma_j(u^*) = 0,\\\quad j =
1,\ldots ,k^2.
\end{multline*}
The function
$$
\sigma = f(s+\omega(x) t,x,\alpha,\beta)\cdot
g(x+iy,\alpha',\beta',\gamma'),
$$
described by  \eqref{section}-\eqref{sectionend} lies in
$\mathcal{L}_4$ for all $\alpha$, $\beta$, $\alpha'$, $\beta'$,
$\gamma'$.
Hence, it expands in terms of the basis $\sigma_j$, and
\begin{equation}
\label{s} a\sigma + \frac{b}{2}(\partial_x-i\partial_y)\sigma +
\frac{c}{2}(\partial_s+\frac{\partial_t}{\omega})\sigma + \\+
\frac{d}{2}(\partial_x+i\partial_y)\sigma = 0.
\end{equation} at $u^*$. Put $L = \frac{b}{2}(\partial_x-i\partial_y)
+\frac{c}{2}(\partial_s+\frac{\partial_t}{\omega}) +
\frac{d}{2}(\partial_x+i\partial_y)$ and rewrite \eqref{s} as
\begin{multline}
\label{a1} L\log((\alpha,\alpha')\cdot\theta_{M})(u^*) = -a - L\log
((\beta,\beta')\cdot \theta_{M})(u^*) -
\\
- L\log((\gamma,\gamma')\cdot\theta_{M})(u^*) -
L\log((\delta,\delta')\cdot\theta_{M})(u^*).
\end{multline}
Here $((\lambda,\mu)\cdot\theta_M)$ is the action described by
\eqref{preprod}.
For all $u,\alpha,\alpha'$ there exist $\beta$, $\beta'$, $\gamma'$,
such that
\begin{equation}
\label{a2} ((\beta,\beta')\cdot
\theta_{M})(u)\times((\gamma,\gamma')\cdot
\theta_{M})(u)\times((\delta,\delta')\cdot \theta_{M})(u)\neq 0.
\end{equation}
It follows from \eqref{a1} and \eqref{a2} that
\begin{equation}
\label{xi} \xi(\alpha,\alpha') =
L\log((\alpha,\alpha')\cdot\theta_{M})(u^*)
\end{equation}
is an entire function of $(\alpha,\alpha')\in \CC^2$.
By \eqref{mult1}-\eqref{mult5} the function $\xi(\alpha,\alpha')$
satisfies the following periodicity conditions:
\begin{eqnarray}
\label{localperiod}
\xi(\alpha + 1,\alpha') & = & \xi(\alpha,\alpha'),\\
\label{localperiod2}
\xi(\alpha + \omega(x^*),\alpha') & = & \xi(\alpha,\alpha') - 2\pi i c - \pi i\frac{(b+d)}{2}\omega(x^*),\\
%\end{eqnarray}
%\begin{eqnarray}
\label{localperiod3}
\xi(\alpha,\alpha' + 1) & = & \xi(\alpha,\alpha'),\\
\label{localperiod4} \xi(\alpha ,\alpha' + i) & = &
\xi(\alpha,\alpha') - 2\pi i b.
\end{eqnarray}
Therefore, derivatives  $\partial_{\alpha}\xi$ and
$\partial_{\alpha'}\xi$ are doubly periodic entire functions. This
means that they are constant and $\xi = A\alpha + B\alpha' + C$.
From \eqref{localperiod} and \eqref{localperiod3} it follows that $A
= B = 0$, and  $\xi\equiv C$. From
\eqref{localperiod2},\eqref{localperiod4} it follows that
$$
b=2\pi i c + \pi i\frac{(b+d)}{2}\omega(x^*) = 0.
$$
Then
\begin{multline}
\label{d} \xi = C =
c\left[\frac{(\partial_s+\frac{\de_t}{\omega}\theta)(s+\omega(x)t+\alpha,\omega(x))}
{\theta(s+\omega(x)t+\alpha,\omega(x))}\right]_{u=u^*} +
\\
\frac{d}{2}\left[\frac{(\partial_x\theta)(s+\omega(x)t+\alpha,\omega(x))}
{\theta(s+\omega(x)t+\alpha,\omega(x))}\right]_{u=u^*} .
\end{multline}
Here we used implicitly the Cauchy--Riemann conditions
$$
(\de_x+i\de_y)\theta(x+iy,i) = 0.
$$
Denote by $D$ the differentiation with respect to $s+\omega t$:
$$
D = \frac{1}{2}\left(\de_s+\frac{\de_t}{\omega}\right) .
$$
It follows from \eqref{classheat} that
\begin{equation}
\label{dd}
\partial_x\theta(s+\omega t,\omega)  = \frac{1}{4\pi i}(D^2\theta)(s+\omega t,\omega)  + t (D\theta)(s+\omega t,\omega).
\end{equation}
Inserting \eqref{dd} into \eqref{d} and taking
$$
(D\theta)(s+\omega t +\alpha,\omega)= \de_{\alpha}\theta(s+\omega t
+ \alpha,\omega),
$$
into account, we find that   ïîëó÷èì, ÷òî  ôóíêöèÿ
$\theta(s^*+\omega t^* + \alpha,\omega)$ as a function of $\alpha$
satisfies a linear ordinary differential equation with constant
coefficients:
$$
\frac{d}{2}\cdot{\chp{\omega(x^*)}{x}}\left(\frac{1}{4\pi i}\theta''
+ t^*\theta'\right)  + c\theta' - C \theta = 0.
$$
Writing down its general solution, we can easily verify that this
leads to a contradiction with the periodicity conditions
\eqref{classperiod}-\eqref{classperiod2} for the theta function;
thereby, $c = d = C = 0$. It follows from  \eqref{s} that $a=0$.

We find that the tuple $a,b,c,d$ of constants is trivial, and the
matrix $\widetilde{J}$ is of maximal rank. Since the point $u^*$ is
chosen arbitrarily, the rank of $\varphi_k$ is equal to $4$
everywhere.  The proof of the theorem is complete.

\section{Embedding is symplectic}

For all $\left\{A,B\right\}$ the total space $M$ of the bundle  is a
symplectic manifold; a symplectic form can be, for instance,
$\omega_M = dx\wedge dy + ds\wedge dt$. In this section we prove the
following statement.

\begin{theorem}
\begin{enumerate}
\item If the mapping $\varphi_k$ is an embedding then it induces a symplectic structure on M.

\item
    The induced symplectic form is cohomologous to $k\cdot
    \omega_M$.
\end{enumerate}
\end{theorem}

{\sc Äîêàçàòåëüñòâî.}\quad In the definition of $\varphi_k$ choose
the functions
$$
\theta^p_k(s+\omega t, \omega)\cdot \theta^{q}_k(x+iy, i);\quad p,q
= 1,\ldots ,k.
$$
as a basis for the space $\mathcal{L}_k$ of theta functions.

Observe that $\varphi_k$ is the composition of the Segre mapping
$\sigma_k:\mathbb{CP}^{k-1} \times \mathbb{CP}^{k-1} \rightarrow
\mathbb{CP}^{k^2-1}$, which is defined in the homogeneous
coordinates as
$$
\sigma_k\left([z^1:\ldots :z^k],[w^1:\ldots :w^k]\right) = [z^1w^1:
z^1w^2: \ldots :z^kw^{k-1}: z^kw^k]
$$ and the map
$\psi_k:M\rightarrow\mathbb{CP}^{k-1}\times\mathbb{CP}^{k-1}$,
$\psi_k = (\psi_k',\psi_k'')$, with
$$
\psi_k'(x,s,t) = [\theta^1_k(s+\omega t,
\omega):\ldots:\theta^k_k(s+\omega t, \omega)],
$$
$$
\psi_k''(x,y) = [\theta^1_k(x+iy,i):\ldots:\theta^k_k(x+iy,i)].
$$
Thus, $\varphi_k = \sigma_k \circ \psi_k$. Denote by  $\Omega'$ the
symplectic form (associated to the Fubini–-Study metric) on the
first factor of $\mathbb{CP}^k~\times~\mathbb{CP}^k$, by $\Omega''$
on the second factor. Then $\Omega'+\Omega''$ is a symplectic form
on the product. Since the Segre mapping is a holomorphic embedding,
it suffices to prove that the induced form $\psi^*_k(\Omega'+
\Omega'')$ is symplectic.

Recall that in Section $3$ we calculated the generators $dx$, $dy$,
$\omega_1$, $\omega_2$ for the algebra of left-invariant forms. The
mapping $\psi_k''$ is a holomorphic embedding of the complex torus
into $\mathbb{CP}^k$ by the classical Lefschetz theorem. Hence,
$$
{(\psi_k'')^*(x,y)} (\Omega'') = \mu\cdot dx\wedge dy,
$$
where $\mu\neq 0$ everywhere on $M$. Put
$$
{(\psi_k')^*(x,s,t)} (\Omega') = f\cdot \omega_1\wedge dx + g\cdot
\omega_2\wedge dx + h\cdot ds\wedge dt
$$
for some functions $f,g,h$ on $M$; this is the general form of a
$2$-form induced by the map $\psi_k'$, which depends on $x,s,t$.

Observe that for every fixed $x$ the map $\psi_k'$ is a holomorphic
embedding as well; thereby,
$$
{(\psi_k')}^* (\Omega') =  \nu \cdot ds\wedge dt,
$$
where $\nu\neq 0$ everywhere on $M$. This implies that $h\equiv
\nu$. Putting everything together, we obtain
\begin{multline*}
\left(\psi^*_k(\Omega'+\Omega'')\right)^2 = \left({(\psi_k')}^*
(\Omega') + {(\psi_k'')}^*
(\Omega'')\right)^2 =\\
 = \left(f\cdot \omega_1\wedge dx + g\cdot \omega_2\wedge dx + \nu\cdot ds\wedge dt + \mu\cdot dx\wedge dy \right)^2.
\end{multline*}
and multiplying out
$$
\left(\psi^*_k(\Omega'+\Omega'')\right)^2 = 2\mu\nu \cdot dx\wedge
dy\wedge ds \wedge dt.
$$
The last equality is equivalent to the non-degeneracy condition for
the induced forms. The closedness follows since the differential
commutes with $\psi_k^*$. Thus, $\psi^*_k(\Omega'+ \Omega'')$ is a
symplectic form. We have proved claim $(1)$ of the theorem.

Let us prove claim $(2)$. Denote by $L$ the bundle given by the
multipliers  \eqref{mult1}-\eqref{mult5}. While establishing the
embedding we observed that the theta functions of degree $k$ are the
sections of $L^{\otimes k}$.

Recall that every complex line bundle over $M$ is induced by the
universal bundle over $\mathbb{CP}^n$ via a mapping of $M$ into the
complex projective space. Consequently, the bundle $L^{\otimes k}$
and its curvature form are the images of the universal bundle and
its curvature form, which is the Fubini–-Study form. Recall also
that the first Chern class of a line bundle is realized precisely by
the curvature form. Hence, the cohomology class of the induced form
coincides with the first Chern class $c_1(L^{\otimes k}) = k\cdot
c_1(L)$ and we must prove that
$$
c_1(L) = [ dx\wedge dy + ds\wedge dt].
$$

Use the \v{C}ech cohomology theory to calculate $c_1(L)$. Cover
$\mathbb{R}^4$ by the open sets
$$
U_{\lambda} = \lambda\cdot U_0,\quad \lambda\in\Gamma.
$$
To this end, dilate the set
$$
U_0 = \{|u^k|<3/4\}.
$$
by the translations of $\Gamma$. Observe that this is a good cover:
all nonempty finite intersections are diffeomorphic to
$\mathbb{R}^4$. Therefore, the cohomologies of the nerve of this
cover are isomorphic to the cohomologies of  $M$.

We can express the transition functions
$g_{\lambda\mu}:U_{\lambda}\cap U_{\mu}\rightarrow \mathbb{C}^*$ in
terms of the multipliers
\begin{equation}
\label{g} g_{\lambda\mu}(u) = e_{\lambda}(u)\cdot
e_{\mu^{-1}}(\mu\cdot u);\quad \lambda,\mu\in\Gamma.
\end{equation}
The nerve $N(\mathcal{U})$ of a minimal subcover of $U_\lambda$ is
homeomorphic to ${M}$, and its cohomology with coefficients in
$\mathbb{Z}$ coincides with $H^*({M};\mathbb{Z})$. The cocycle
$z_{\lambda\mu\nu} \in C^2(\mathcal{U};\mathbb{Z})$
\begin{equation}
\label{cocycle} z_{\lambda\mu\nu} = \frac{1}{2\pi
i}(\log(g_{\lambda\mu}) + \log(g_{\mu\nu}) - \log(g_{\nu\lambda}))
\end{equation}
by definition realizes the first Chern class of the bundle $L$. This
formula defines the value of $z$ on the $2$-dimensional simplex
$(\lambda,\mu,\nu)\in N(\mathcal{U})$.

It is easy to establish that

$(1)$ $H^2(M;\RR) \cong \RR^4$, if  $M$ is the Kodaira-–Thurston
manifold:
$$
A = \matriza{1}{\lambda}{0}{1},\quad B =E,
$$
and the generators for $H^2(M;\RR)$ can be chosen to be $[dx\wedge
dy]$, $[ds\wedge dt]$, $[dy\wedge dt]$, $[(ds-\lambda x dt)\wedge
dx]$;

$(2)$ in all remaining cases $H^2(M;\RR) \cong \RR^2$ and the
generators are $[dx\wedge dy]$, $[ds\wedge dt]$.

Therefore, for all manifolds $M$ tori $T_{ab}$ and $T_{cd}$ formed
by the commuting translations \eqref{gammaA}-\eqref{gammaBB} in
$\Gamma$ are dual cycles to $[dx\wedge dy]$, $[ds\wedge dt]$. In the
case of the Kodaira-–Thurston manifold we obtain two more cycles,
$T_{bd}$ and $T_{ac}$.

Define the functions $f_\lambda(u)$ so as
\begin{equation}
\label{f} e_\lambda(u) = e^{2\pi i f_\lambda(u)}.
\end{equation}
By \eqref{g}--\eqref{f}
\begin{equation*}
%\label{c1}
c_1([T_{\lambda\mu}]) = f_\mu(u) + f_\lambda(\mu \cdot u) -
f_\lambda(u) - f_\mu(\lambda\cdot u).
\end{equation*}
Evaluating the first Chern class on $T_{ab}$ and  $T_{cd}$, we
obtain
\begin{equation}
\label{c1calc} c_1([T_{ab}]) = c_1([T_{cd}]) = 1;
\end{equation}
In the case of the Kodaira–-Thurston manifold,\begin{equation}
\label{c1calc2} c_1([T_{bd}]) =  c_1([T_{ac}]) = 0.
\end{equation}

Since the manifold $M$ is a homogeneous space of a real Lie group,
all elements of  $H^2({M};\mathbb{R})$ can be realized by
left-invariant forms dual to the basis $2$-cycles. From
\eqref{c1calc} and \eqref{c1calc2} it follows that $c_1(L) =
[dx\wedge dy + ds\wedge dt]$. The proof of the theorem is complete.

\end{document}